\documentclass[11pt,letterpaper]{amsart}
\usepackage[usenames,dvipsnames]{color}
\usepackage{amsthm,amsfonts,amssymb,amsmath,amsxtra}
\usepackage[all]{xy}
\SelectTips{cm}{}
\usepackage{xr-hyper}
\usepackage[pdfpagelabels,pdftex,hidelinks]{hyperref}
\hypersetup{
  colorlinks   = true, 
  urlcolor     = blue, 
  linkcolor    = blue, 
  citecolor   = blue, 
pdftitle={},
  pdfauthor={},
  pdfsubject={},
  pdfkeywords={},
  breaklinks=true,
  bookmarksopen=true,
  bookmarksnumbered=true,
  pdfpagemode=UseOutlines,
  plainpages=false}
\usepackage{cleveref}

\usepackage{verbatim}
\usepackage{caption}

\usepackage{lineno}

\usepackage{mathrsfs}

\RequirePackage{xspace}
\RequirePackage{etoolbox}
\RequirePackage{varwidth}
\RequirePackage{enumitem}
\RequirePackage{mathtools}
\RequirePackage{longtable}
\RequirePackage{multirow}

\def\ge{\geqslant}
\def\le{\leqslant}
\def\a{\alpha}
\def\b{\beta}

\def\G{\Gamma}
\def\d{\delta}
\def\D{\Delta}

\def\t{\tau}
\def\th{\theta}

\def\l{\lambda}
\def\z{\zeta}

\def\i{^{-1}}

\def\<{\langle}
\def\>{\rangle}

\newcommand{\fkX}{\ensuremath{\mathfrak{X}}\xspace}

\newcommand{\BA}{\ensuremath{\mathbb {A}}\xspace}

\newcommand{\BF}{\ensuremath{\mathbb {F}}\xspace}
\newcommand{{\BG}}{\ensuremath{\mathbb {G}}\xspace}

\newcommand{{\BK}}{\ensuremath{\mathbb {K}}\xspace}

\newcommand{\BQ}{\ensuremath{\mathbb {Q}}\xspace}
\newcommand{\BR}{\ensuremath{\mathbb {R}}\xspace}

\newcommand{\BZ}{\ensuremath{\mathbb {Z}}\xspace}

\newcommand{\CA}{\ensuremath{\mathcal {A}}\xspace}
\newcommand{\CB}{\ensuremath{\mathcal {B}}\xspace}

\newcommand{\CE}{\ensuremath{\mathcal {E}}\xspace}
\newcommand{\CF}{\ensuremath{\mathcal {F}}\xspace}
\newcommand{\CG}{\ensuremath{\mathcal {G}}\xspace}

\newcommand{\CI}{\ensuremath{\mathcal {I}}\xspace}

\newcommand{\CK}{\ensuremath{\mathcal {K}}\xspace}
\newcommand{\CL}{\ensuremath{\mathcal {L}}\xspace}
\newcommand{\CM}{\ensuremath{\mathcal {M}}\xspace}

\newcommand{\CO}{\ensuremath{\mathcal {O}}\xspace}
\newcommand{\CP}{\ensuremath{\mathcal {P}}\xspace}

\newcommand{\CR}{\ensuremath{\mathcal {R}}\xspace}

\newcommand{\CT}{\ensuremath{\mathcal {T}}\xspace}

\newcommand{\GL}{\mathrm{GL}}

\DeclareMathOperator{\Hom}{Hom}

\let\Im\relax
\DeclareMathOperator{\Im}{Im}

\newcommand{\ov}{\overline}

\def\tW{\widetilde W}

\def\brk{{\breve k}}

\def\COk{{\CO_{\brk}}}

\def\pr{{\rm pr}}
\def\tPhi{\widetilde \Phi}

\def\ind{{\rm ind}}

\def\bx{{\mathbf x}}

\def\der{{\rm der}}
\def\ov{\overline}

\def\hSigma{{\hat \Sigma}}

\def\hXi{{\widehat \Xi}}

\def\re{{\rm re}}
\def\dz{{\dot z}}

\def\hCI{{\widehat{\CI}}}
\def\hSigma{{\widehat \Sigma}}

\newtheorem{theorem}{Theorem}
\newtheorem{proposition}[theorem]{Proposition}
\newtheorem{lemma}[theorem]{Lemma}

\newtheorem{corollary}[theorem]{Corollary}

\theoremstyle{definition}

\numberwithin{equation}{section}
\numberwithin{theorem}{section}

\setitemize[0]{leftmargin=*,itemsep=\the\smallskipamount}
\setenumerate[0]{leftmargin=*,itemsep=\the\smallskipamount}

\renewcommand{\to}{%
   \ifbool{@display}{\longrightarrow}{\rightarrow}%
   }
\let\shortmapsto\mapsto
\renewcommand{\mapsto}{%
   \ifbool{@display}{\longmapsto}{\shortmapsto}%
   }
\newlength{\olen}
\newlength{\ulen}
\newlength{\xlen}
\newcommand{\xra}[2][]{%
   \ifbool{@display}%
      {\settowidth{\olen}{$\overset{#2}{\longrightarrow}$}%
       \settowidth{\ulen}{$\underset{#1}{\longrightarrow}$}%
       \settowidth{\xlen}{$\xrightarrow[#1]{#2}$}%
       \ifdimgreater{\olen}{\xlen}%
          {\underset{#1}{\overset{#2}{\longrightarrow}}}%
          {\ifdimgreater{\ulen}{\xlen}%
             {\underset{#1}{\overset{#2}{\longrightarrow}}}
             {\xrightarrow[#1]{#2}}}}%
      {\xrightarrow[#1]{#2}}
   }
\makeatother
\newcommand{\xyra}[2][]{%
   \settowidth{\xlen}{$\xrightarrow[#1]{#2}$}%
   \ifbool{@display}%
      {\settowidth{\olen}{$\overset{#2}{\longrightarrow}$}%
       \settowidth{\ulen}{$\underset{#1}{\longrightarrow}$}%
       \ifdimgreater{\olen}{\xlen}%
          {\mathrel{\xymatrix@M=.12ex@C=3.2ex{\ar[r]^-{#2}_-{#1} &}}}%
          {\ifdimgreater{\ulen}{\xlen}%
             {\mathrel{\xymatrix@M=.12ex@C=3.2ex{\ar[r]^-{#2}_-{#1} &}}}
             {\mathrel{\xymatrix@M=.12ex@C=\the\xlen{\ar[r]^-{#2}_-{#1} &}}}}}%
      {\mathrel{\xymatrix@M=.12ex@C=\the\xlen{\ar[r]^-{#2}_-{#1} &}}}%
   }
\makeatletter
\newcommand{\xla}[2][]{%
   \ifbool{@display}%
      {\settowidth{\olen}{$\overset{#2}{\longleftarrow}$}%
       \settowidth{\ulen}{$\underset{#1}{\longleftarrow}$}%
       \settowidth{\xlen}{$\xleftarrow[#1]{#2}$}%
       \ifdimgreater{\olen}{\xlen}%
          {\underset{#1}{\overset{#2}{\longleftarrow}}}%
          {\ifdimgreater{\ulen}{\xlen}%
             {\underset{#1}{\overset{#2}{\longleftarrow}}}
             {\xleftarrow[#1]{#2}}}}%
      {\xleftarrow[#1]{#2}}
   }
\newcommand{\isoarrow}{%
   \ifbool{@display}{\overset{\sim}{\longrightarrow}}{\xrightarrow\sim}%
   }

\newcommand{\sm}{{\,\smallsetminus\,}}

\usepackage{upgreek}
\makeatletter
\newcommand{\colim@}[2]{%
  \vtop{\m@th\ialign{##\cr
    \hfil$#1\operator@font lim$\hfil\cr
    \noalign{\nointerlineskip\kern1.5\ex@}#2\cr
    \noalign{\nointerlineskip\kern-\ex@}\cr}}%
}
\newcommand{\colim}{%
  \mathop{\mathpalette\colim@{\rightarrowfill@\textstyle}}\nmlimits@
}
\newcommand{\prolim@}[2]{%
  \vtop{\m@th\ialign{##\cr
    \hfil$#1\operator@font lim$\hfil\cr
    \noalign{\nointerlineskip\kern1.5\ex@}#2\cr
    \noalign{\nointerlineskip\kern-\ex@}\cr}}%
}
\newcommand{\prolim}{%
  \mathop{\mathpalette\colim@{\leftarrowfill@\textstyle}}\nmlimits@
}
\makeatother

\DeclareMathOperator{\Perf}{Perf}

\begin{document}

\title[]{Deep level Deligne-Lusztig representations of Convex type} 

\author[Ben Liu]{Ben Liu}
\address{Academy of Mathematics and Systems Science, Chinese Academy of Sciences, Beijing 100190, China}
\email{liubenmath@gmail.com}

\author[Sian Nie]{Sian Nie}
\address{Academy of Mathematics and Systems Science, Chinese Academy of Sciences, Beijing 100190, China}

\address{School of Mathematical Sciences, University of Chinese Academy of Sciences, Chinese Academy of Sciences, Beijing 100049, China}
\email{niesian@amss.ac.cn}

\author{Zhihang Yu}
\address{The Second Affiliated Hospital of Chongqing Medical University, Chongqing 400010, China}
\email{yuzhihang@amss.ac.cn}

\begin{abstract}
In this paper, we establish a dimension formula for intertwining spaces of Lusztig's deep level Deligne-Lusztig representations attached to convex Borel subgroups. This is achieved by revising the previous construction of these representations by the authors. As an application, we provide a construction of a large class of irreducible supercuspidal representations of unramified $p$-adic groups for all primes $p$, in a purely geometric way.  
\end{abstract}
\maketitle

\section{Introduction}

Let $G$ be a connected reductive group defined over a non-archimedean local field $k$ with ring of integers $\CO_k$ and residue characteristic $p$. Fix a point $\bx$ in the Bruhat-Tits building and let $\CG_{\bx, 0}$ be the corresponding connected parahoric group scheme. Let $T \subseteq G$ be an unramified elliptic maximal torus over $k$ whose apartment contains $\bx$. Let $Z$ denote the center of $G$.

As in the classical Deligne-Lusztig theory \cite{DL76}, to each smooth character $\phi$ of $T(k)$ and a Borel subgroup $B \supseteq T$ one can attach a virtual smooth $Z(k) \CG_{\bx, 0}(\CO_k)$-module \[R(\phi) = R_{T \subseteq B} (\phi)\] using the cohomological induction on a suitable subvariety of the truncated jet scheme of $\CG_{\bx, 0}$. These virtual modules are referred to as deep level Deligne-Lusztig representations. They were first introduced by Lusztig \cite{Lu79} to provide a geometric realization of G\'{e}rardin's construction of supercuspidal representations in \cite{Ger75a}. Over the past decades, the study of $R(\phi)$ has undergone rapid progress and shown striking applications to the local Langlands program, see \cite{BW16}, \cite{CO25a}, \cite{BC24} and \cite{INY25}. However, almost all known results on $R(\phi)$ are based upon a key assumption that $\phi$ admits a Howe factorization in the sense of \cite{Kal}, which may fail when $p$ is a bad prime for $G$ or $p$ divides $|\pi_1(G_\der)|$.

\subsection{The main results}
In this paper, we fulfill Lusztig's geometric approach to supercuspidals for all smooth characters $\phi$. Our first main result is a dimension formula of virtual intertwining spaces.  
\begin{theorem} \label{intro-inter}
    Assume the Borel subgroup $B \supseteq T$ is convex. For $w \in G(k)$ we have 
    \[\dim I_w(R(\phi)) = |\{z \in (W_\bx z_w W_\bx)^F; {}^z \phi |_{\CT(\CO_k)} = \phi|_{\CT(\CO_k)}\}|,\] where $W_\bx$ is the absolute Weyl group of the $\CO_k$-model $\CT$ of $T$ in $\CG_{\bx, 0}$, $F$ is the Frobenius map of $G$, and $z_w$ is the Iwahori-Weyl group element associated to $w$.
\end{theorem}
The notion of convex Borel subgroups, or Weyl group elements, was first introduced in \cite{IN25} to compute Fargues--Scholze parameters of modular supercuspidal representations. It is proved in \cite{NTY25} that convex Borel subgroups in Theorem \ref{intro-inter} always exist. Note that the Borel subgroups of Coxeter type considered in \cite{Chan20}, \cite{CI23}, \cite{DI24} and \cite{IN24} are special cases of convex Borel subgroups.   

Under a mild triviality condition of $\phi$ (see \cite[Remark 1.1]{LNY26}), such a dimension formula was first established in \cite[Theorem 1.2]{LNY26}, by using an alternative realization $\CR(\phi)$ of $R(\phi)$ introduced in \cite{CS17}, \cite{CS23} and \cite{Nie24}. 

As a special case, we deduce the Mackey inner product formula for $R(\phi)$ associated to convex Borels. 
\begin{corollary} \label{intro-Mackey}
    Assume $B$ is convex. Then \[\dim \Hom_{\CG_{\bx, 0}(\CO_k)}(R(\phi), R(\phi)) = |\{z \in W_\bx^F; {}^z \phi |_{\CT(\CO_k)} = \phi|_{\CT(\CO_k)}\}|.\] In particular,
    $\pm R(\phi)$ is irreducible if and only if the stabilizer of $\phi|_{\CT(\CO_k)}$ in $W_\bx^F$ is trivial.
\end{corollary}
When $\phi$ admits a Howe factorization, this Mackey formula is proved in \cite{Chan24} for all Borel subgroups containing $T$.

As in \cite{LNY26}, we obtain a purely geometric construction of a large class of irreducible supercuspidal representations.
\begin{theorem} [cf. {\cite[Theorem 1.4]{LNY26}}] \label{intro-sum}
    Assume $\pm R(\phi)$ is a genuine representation. Then \[\pi(\phi) := \pm \text{c-} \ind_{Z(k) \CG_{\bx, 0}(\CO_k)}^{G(k)} R(\phi)\] is a direct sum of irreducible supercuspidals.
\end{theorem}

\begin{theorem} [cf. {\cite[Theorem 1.5]{LNY26}}] \label{intro-irr-sup}
    Assume  the stabilizer of $\phi|_{\CT(\CO_k)}$ in $\pr(\tW_{[\bx]}^F)$ is trivial. Then \[\text{$\pi(\phi)$ is an irreducible supercuspidal.}\] Here $[\bx]$ is the image of $\bx$ in the reduced building, $\tW_{[\bx]}$ is the stabilizer of $[\bx]$ in the Iwahori-Weyl group and $\pr(\tW_{[\bx]})$ denotes its natural image in the Weyl group.
\end{theorem}
Similar results are obtained in \cite{Chan20}, \cite{CS17}, \cite{CI23}, \cite{CO25a}, \cite{Nie24} and \cite{INT24}, under various mild conditions on $\phi$ or the residue field of $k$.

\subsection{The strategy and organization}
The proof of Theorem \ref{intro-inter} consists of two ingredients. The first one is a revised construction of deep level Deligne-Lusztig representations, see Proposition \ref{new-constr}. It extends the construction of \cite{Nie24} and \cite{LNY26} in a natural way, and makes sense for all smooth characters $\phi$. The second ingredient is the notion of convex Borel subgroups introduced in \cite{IN25}. It enables us to show that Lusztig's construction and the new construction are equivalent in the sense of cohomological complexes up to shifts. This is achieved by adapting methods and techniques developed in \cite{IN25} and \cite{LNY26}.

The paper is organized as follows. In \S\ref{sec:pre} we collect basic notions and introduce the notion of convex Borel subgroups. In \S\ref{sec:dl-DL} we introduce the notion of deep level Deligne-Lusztig representations in the sense of Lusztig. In \S\ref{sec:new-constr}, we show that in the convex case Lusztig's construction coincides with our new construction up to cohomological shifts. In \S\ref{sec:inter}, \S\ref{sec:Sigma} and \S\ref{sec:dim-inter} we establish the dimension formula of intertwining spaces in Theorem \ref{intro-inter}. In \S\ref{sec:sum} we establish a vanishing result for deep Deligne-Lusztig representations, which leads to Theorem \ref{intro-sum}.

\section{Preliminaries} \label{sec:pre}
Let $k$ be a local non-Archimedean field with residue field $\BF_q$ of characteristic $p$. We write $\brk$ for the completion of a maximal unramified extension of $k$, and $F$ for the Frobenius automorphism of $\brk$ over $k$. We denote by $\CO_k$ and $\COk$ the integer rings of $k$ and $\brk$ respectively. Let $\varpi$ be a fixed uniformizer of $k$. We also fix a prime $\ell\neq p$.

\subsection{Positive loops}
Let $\Perf_{\ov\BF_q}$ denote the category of perfect $\ov\BF_q$-schemes. Let $\fkX$ be an $\COk$-scheme. The positive loop functor (also called jet scheme) of $\fkX$ is defined by
\[
L^+ \fkX \colon \Perf_{\ov\BF_q}^{\rm op} \to {\rm Sets}, \quad R \mapsto \fkX(W_{\COk}(R)),
\]
where $W_{\COk}(R) = R[[\varpi]]$ if ${\rm char}\, \brk > 0$ and $W_{\COk}(R) = \COk \otimes_{\breve \BZ_p} W(R)$ otherwise. Here $W(R)$ denotes the $p$-typical Witt-vectors of $R$. In this paper, we only consider schemes $\fkX$ which are affine and of finite type over $\COk$, in which case $L^+\fkX$ is representable by an affine perfect scheme over $\ov\BF_q$.

\subsection{Moy-Prasad quotients}
Let $G$ be a connected reductive group defined over $k$. Let $\bx$ be a point in the Bruhat-Tits building $\CB(G, k)$ of $G$ over $k$. Let $\CT$ be the corresponding integral model for $T$. Bruhat--Tits theory associates to each $0 \le r \in \widetilde \BR := \BR \sqcup \{s+; s \in \BR\}$ a Moy-Prasad $\CO_k$-group model of $G$. For any $0 \le s \le r \in \widetilde \BR$, $L^+(\CG_{\bx,r+} \times_{\CO_k} \CO_{\breve k}) \subseteq L^+(\CG_{\bx,s} \times_{\CO_k} \CO_{\breve k})$ is a subgroup, and we put
\[G_{s:r} = L^+(\CG_{\bx,s} \times_{\CO_k} \CO_{\breve k}) / L^+(\CG_{\bx,r+} \times_{\CO_k} \CO_{\breve k}),\]
which is a pfp perfectly smooth affine group scheme over $\ov\BF_q$. As the Frobenius automorphism $F$ fixes $\bx$, it acts naturally on $G_{s:r}$ and hence endows it with an $\BF_q$-rational structure. In particular, $G_{s:r}(\BF_q) = \CG_{\bx,s}(\CO_k)/\CG_{\bx,r+}(\CO_k)$.

Let $H \subseteq G$ be a closed subgroup. Following \cite[\S2.5]{CI21} one may construct an $\ov\BF_q$-rational closed subgroup $H_{s:r} \subseteq G_{s:r}$ in a similar way. We put $H_r = H_{0:r}$ for simplicity. If, moreover, $H$ is a $k$-rational subgroup, then $H_{s:r}$ is defined over $\BF_q$, and we still denote by $F$ the induced Frobenius automorphisms on $H$ or $H_{s:r}$. We put $H = H(\brk)$ and $H_{s:r} = H_{s:r}(\ov\BF_q)$ for simplicity.

\subsection{The maximal torus}
Let $T \subseteq G$ be a $k$-rational and $\brk$-split elliptic maximal torus such that $\bx \in \CA(T, \brk) \cap \CB(G, k)$. Here $\CA(T, \brk)$ denotes the apartment of $T$ over $\brk$. By fixing a special point $o \in \CA(T, \brk)$ fixed by the Frobenius of the unique quasi-split inner form of $G$, we identify $\CA(T, \brk)$ with the underlying affine space of the real vector space $X_*(T) \otimes \BR$. By this identification, the action of $F$ on $\CA(T, \brk)$ induces an affine transformation of $X_*(T) \otimes \BR$, which we still denote by $F$. 

Let $N_T$ be the normalizer of $T$ in $G$. We have the finite Weyl group $W = N_T(\brk) / T(\brk)$ and the Iwahori-Weyl group $\tW = N_T(\brk) / T(\COk)$. As the Frobenius map $F$ preserves $T$, it acts on $W$ and $\tW$ in a natural way. There is a natural action of $\tW \rtimes \<F\>$ on $X_*(T) \otimes \BR$ by affine transformations.

Let $\Phi = \Phi(G_\brk, T_\brk)$ denote the root system of $T_\brk$ in $G_\brk$. For $\a \in \Phi$ we write $G^\a \subseteq G$ and $\a^\vee: \BG_m \to T_\brk$ for the corresponding root subgroup and coroot respectively. We fix a normalized Chevalley system for $(G, T)$ compatible with the point $o$, which is a collection of isomorphisms $u_\a : \BG_a \to G^\a$ for $\a \in \Phi$ subject to certain properties, see \cite{BT2} or \cite{Ad}. 

Let $\tPhi = \tPhi_\re \sqcup\, \BZ$ be the set of affine roots, where $\tPhi_\re = \Phi \times \BZ$ denotes the set of real affine roots. Let $f \in \tPhi$. We write $\a_f \in \Phi \sqcup \{0\}$, $n_f \in \BZ$ such that $f = (\a_f,n_f) = \a_f + n_f$. We view $f$ as an affine function on $X_*(T) \otimes \BR$ given by $f(v) = \a_f(v) + n_f$. Then the action of $\tW \rtimes \<F\>$ on $X_*(T) \otimes \BR$ induces an action on $\tPhi$ by $x(f)(v) = f(x\i(v))$.

We fix a basis $(\l_i)_{i=1}^l$ of the cocharacter group $X_*(T)$ of $T$. If $f \in \BZ_{\ge 1}$, we put $\CG^f =\CT(1 + \varpi^f\CO_\brk)$ and define \[u_f: \BA^l \to \CG^f, \quad (x_i)_i \mapsto \prod_i \l_i(1 + \varpi^{n_f} [x_i]),\] where $[x_i] = x_i$ in the equal characteristic case and $[x_i]$ is the Teichm\"{u}ller lift of $x_i \in \ov\BF_q$ in the mixed characteristic case. If $f \in \tPhi_\re$, we put $\CG^f = u_{\a_f}(\varpi^{n_f} \CO_\brk)$ and define \[u_f: \BA^1 \to \CG^f, \quad x \mapsto u_{\a_f}(\varpi^{n_f} [x]).\] We put $G^f = \Im (u_f)$ for $f \in \tPhi_\re \sqcup \BZ_{\ge 1}$ with $f(\bx)\ge0$. The parahoric subgroup $G_{\bx,0} := \CG_{\bx, 0}(\COk)$ attached to $\bx$ is generated by $\CT(\COk)$ and $\CG^f$ for $f \in \tPhi_\re$ such that $f(\bx) \ge 0$. Moreover, we have the Bruhat decomposition \[G = G(\brk) = \sqcup_{z \in W_\bx \backslash \tW / W_\bx} G_{\bx, 0} z G_{\bx, 0},\] where $W_\bx \subseteq \tW$ is the Weyl group of $G_{\bx, 0}$.

\medskip

Let $B = T U \subseteq G$ be a Borel subgroup with unipotent radical $U$. We denote by $\Phi^+$ the set of (positive) roots which appear in $B$. Let $\Pi = \{\a_1, \dots, \a_d\}$ be the set of simple roots of $\Phi^+$. For $\nu \in \BZ\Phi$ we write $\nu \prec 0$ if $\nu = c_1 \a_1 + \cdots + c_d \a_d \in \BZ\Phi$ such that there exists $1 \le i \le d$ such that $c_1 = \cdots = c_{i-1} = 0$ and $c_i < 0$. This defines a linear order $\preceq$ on $\BZ\Phi$ by \[\l \prec \l' \Leftrightarrow (\l-\l') \prec 0.\] Let $\leq$ be a linear order on $(\BZ\Phi) \times \BR$ such that $\l + n < \l' + n'$ if and only if either $\l(\bx) + n < \l'(\bx) + n'$ or $\l(\bx) + n = \l'(\bx) + n'$ and $\l \prec \l'$. Let $\tPhi^{>0}$ be the set of affine roots $f > 0$, that is, either $f(\bx) > 0$ or $f(\bx) = 0$ and $\a_f \in \Phi^+$.

Let $r \in \BZ_{\ge 0}$, $\tPhi_{\re,r}=\{f\in\tPhi_\re;0\le f(\bx)\le r\}$, and $\tPhi_{\re, r}^{>0} = \{f \in \tPhi_\re^{>0}; f(\bx) \le r\}$. We denote by $\CG_r^f$ and $G_r^f$ the natural images of $\CG^f$ and $G^f$ in $G_r$ for $f \in \tPhi_{\re,r}$ respectively. For $E \subseteq \tPhi_{\re, r}^{>0} $, we put $G_r^E = \prod_{f \in E}G_r^f$, where the product is taken with respect to the linear order $\le$ on $(\BZ\Phi) \times \BR$.  For $v \in G_r^E$ there are unique elements $v_f \in G_r^f$ such that $v = \prod_{f \in E} v_f$. By abuse of notation, we also view $v_f$ as elements in $\ov\BF_q$ via the bijections $\BA^1 \overset {u_f} \to G_r^f$.

\subsection{Convex Borel subgroups}
Recall that $B = T U$ is a Borel subgroup over $\brk$ and $\Phi^+$ is the set of positive roots associated to $B$.

Let $\a \in \Phi^+$. We define \begin{align*} n_F(\a) &= \min\{i \in \BZ_{\ge 1}; F^i(\a) \in -\Phi^+\}; \\ n_{F\i}(\a) &= \min\{i \in \BZ_{\ge 1}; F^{-i}(\a) \in -\Phi^+\}, \end{align*} which are well-defined since $T$ is elliptic.

We say $B$ is convex if $n_F(\a+\b) \le \max\{n_F(\a), n_F(\b)\}$ and $n_{F\i}(\a+\b) \le \max\{n_{F\i}(\a), n_{F\i}(\b)\}$ for $\a, \b \in \Phi^+$ such that $\a + \b \in \Phi$. By \cite{NTY25}, there always exists a convex Borel subgroup over $\brk$ containing $T$.

\section{Deep level Deligne-Lusztig representations} \label{sec:dl-DL}
Let $r \in \BZ_{\ge 0}$ and $L: G_r \to G_r$ be the Lang map given by $g \mapsto g\i F(g)$. Fix a smooth character $\phi: T^F \to \ov\BQ_\ell^\times$ of depth $r_\phi \le r$.

Let $\a \in \Phi$. We define \[r_\a = 2 s_\a = \min\{i \in \BZ_{\ge 0}; \phi \circ N_F^{F^N} \circ \a^\vee(1 + \varpi^{i+1}\CO_\brk^{F^N}) = \{1\}\}.\] Here $N$ is any positive integer such that $F^N$ acts trivially on $X_*(T)$ and $N_F^{F^N}: T \to T$ is given by $x \mapsto x F(x) \cdots F^{N-1}(x)$. By definition, we have $r_\a = r_{-\a} = r_{F(\a)}$.

Let $f \in \tPhi$. We define $f^\dag = r_f - f$. Here $r_f = 2 s_f = r_{\a_f}$ for $f \in \tPhi_\re$ and $r_f = 2 s_f = -\infty$ if $f \in \BZ$. Let $\CK_{\phi, r}$ (resp. $\CK_{\phi, r}^+$) be the subgroup of $G_r$ generated by $T_r$ (resp. $T_{0+:r}$) and $G_r^f$ for $f \in \tPhi_\re$ such that $f(\bx) \ge s_f$ (resp. $f(\bx) > s_f$).

We recall the following result from \cite[Corollary 3.3]{LNY26}.
\begin{lemma} \label{comm}
    Let $f, f' \in \tPhi_{\re, r}$ with $\a_f + \a_{f'} \neq 0$. Then $[\CG_r^f, \CG_r^{f'}]$ is contained in the product of subgroups $\CG_r^{f''}$ with $f'' = if + jf' + m_{f,f',i,j} \in \tPhi_\re$ with $i, j \in \BZ_{\ge 1}$ and $m_{f,f',i,j} \in \BZ_{\ge 0}$ such that $r_{f''} - m_{f,f',i,j} \le \max\{s_{f''}, r_f, r_{f'}\}$.
\end{lemma}

\subsection{The subgroup $\CI_{\phi, r}^+$}  
For $\a \in \Phi$ with $n \ge 1$  we set \[ \l_{\a^\vee, n}: \BA^1 \to T_r, \quad x \mapsto \a^\vee(1 + \varpi^n [x]).\]

Let $E_\phi = \{f \in \tPhi_\re; f > s_f\}$. Write $E_{\phi, r} = E_\phi \cap \tPhi_{\re, r}^{>0},
\CI_{\phi,r}^+=T_{0+:r}G_r^{E_{\phi,r}},D_{\phi,r}=\{(\a,n)\in\Phi\times\BZ;r_\a+1\le n\le r\}$ and we have the parameterization
\begin{gather*}
\hCI_{\phi,r}^\dag:\BA^{D_{\phi,r}}\times\BA^{E_{\phi,r}}\longrightarrow\CI_{\phi,r}^+,\\
(z_\z,x_f)_{\z,f}\longmapsto
\prod_{\z=(\a,n)\in D_{\phi,r}}\l_{\a^\vee,n}(z_\z)
\prod_{f\in E_{\phi,r}}u_f(x_f).
\end{gather*}
Let $\th_1,\th_2: X \to \CI_{\phi, r}^+$ be any two morphisms. Write \[\th_1 \equiv \th_2 \mod \hCI_{\phi, r}^\dag \quad\text{if }\th_1\i\th_2\text{ factors through }\hCI_{\phi,r}^\dag.\]

Note that $T_0 = T_r / T_{0+:r}$. By Teichm\"{u}ller lift, we also view $T_0$ as the reductive part of $T_r$ and hence $T_r \cong T_0 \times T_{0+:r}$. We put $\phi_0 = \phi|_{T_0^F}$. Let $\CL_\phi$ (resp. $\CL$) be the multiplicative local system over $T_r$ (resp. $T_{0+:r}$) associated to the character $\phi|_{T_r^F}$ (resp. $\phi |_{T_{0+:r}^F}$). 

Let $\t_r: \CI_{\phi, r}^+ \to T_{0+:r}$ be the projection to the first factor. We put $\CF = \t_r^*\CL$.

\begin{lemma} \label{L_phi}
    Let $\a \in \Phi$ and $n \in \BZ_{\ge 1}$. Then $\l_{\a^\vee, n}^* \CL$ is a nontrivial Artin-Schreier sheaf if $n = r_\a > 0$, and $\l_{\a^\vee, n}^* \CL$ is trivial if $n \ge r_\a + 1$. In particular, $\th_1^*\CF\cong \th_2^*\CF$ if $\th_1 \equiv \th_2 \mod \hCI_{\phi, r}^\dag$.
\end{lemma}

\begin{lemma} \label{local-CF}
    The local system $\CF$ on $\CI_{\phi, r}^+$ is multiplicative and $T_0$-equivariant by conjugation.
\end{lemma}
\begin{proof}
    Let $m_1: \CI_{\phi, r}^+ \times \CI_{\phi, r}^+ \to \CI_{\phi, r}^+$ and $m_2: T_{0+:r} \times T_{0+:r} \to T_{0+:r}$ be the multiplication morphisms. By the proof of \cite[Proposition 4.1]{LNY26} we have \[\t_r \circ m_1 \equiv m_2 \circ (\t_r \times \t_r) \mod \hCI_{\phi, r}^\dag.\] Note that $\CL$ is multiplicative. Applying Lemma \ref{L_phi} we have \[m_1^*\CF =m_1^* \t_r^*\CL \cong (\t_r \times \t_r)^* {m_2}^* \CL \cong (\t_r \times \t_r)^* (\CL \boxtimes \CL) \cong \CF \boxtimes \CF.\] Hence $\CF$ is multiplicative.
\end{proof}

For $c \in \ov\BF_q$ let $\mu_c: \BA^1 \to \BA^1$ be the scalar map $z \mapsto c z$.
\begin{lemma} \label{commutator}
    Let $f \in E_\phi \sqcup \BZ_{\ge1} \subseteq \tPhi$ and $y = \prod_{h \in \tPhi_{\re, r}^{>0}} u_h(y_h)$ with $y_h \in \ov\BF_q$ such that \[\text{ $f \ge h^\dag$ if $h < s_h$ and $y_h \neq 0$.}\] Then, for any $\t \in T_{0+:r}$, we have
    \[[-, \t y]\circ u_f \equiv \begin{cases}\l_{\a_f^\vee, r_f} \circ \mu_{y_{f^\dag}}&r_f>0,\\1&r_f=0,\end{cases} \mod \hCI_{\phi, r}^\dag \] 
\end{lemma}
\begin{proof}
    It follows in a similar way as \cite[Corollary 4.2]{LNY26}.
\end{proof}

\subsection{The representation $R_{T_r}^{G_r}(\phi)$}
Let $X \subseteq G_r$ be a locally closed subset normalized by $T_r$. The preimage \[L\i(X)\] admits an action of $G_r^F \times T_r^F$ given by $(g, t): x \mapsto g x t$. Then the $\phi$-isotypic subspace \[H_c^i(L\i(X), \ov\BQ_\ell)[\phi] \subseteq H_c^i(L\i(X), \ov\BQ_\ell)\] is a $G_r^F$-module. We obtain a virtual $G_r^F$-module \[H_c^*(L\i(X), \ov\BQ_\ell)[\phi] := \sum_{i \in \BZ} (-1)^i H_c^i(L\i(X), \ov\BQ_\ell)[\phi].\] In particular, if $X = FU_r$, the virtual $G_r^F$-module \[R_{T_r}^{G_r}(\phi) = H_c^*(L\i(FU_r), \ov\BQ_\ell)[\phi]\] is called a deep level Deligne-Lusztig representation of $G_r^F$. By inflation, we also view $R_{T_r}^{G_r}(\phi)$ as a virtual smooth representation of $G_{\bx, 0}^F$.

Let $Z$ be the center of $G$. As $T$ is elliptic, we have $T^F =Z^F T(\COk)^F$. Thus we can extend each $G_{\bx, 0}^F$-module $H_c^i(L\i(FU_r), \ov\BQ_\ell)[\phi]$ to a $Z^F G_{\bx, 0}^F$-module on which $Z^F$ acts by the character $\phi |_{Z^F}$. In this way, we view $R_{T_r}^{G_r}(\phi)$ as a virtual $Z^F G_{\bx, 0}^F$-module.

\section{Reduction to the locus $L\i(F\CI_{\phi, r}^+)$} \label{sec:new-constr}

In this section, we reduce the computation of $R\G_c(L\i(FU_r), \ov\BQ_\ell)[\phi]$ to that of $R\G_c(L\i(F\CI_{\phi, r}^+), (F\i L)^*\CF)[\phi_0]$, under the assumption that $B$ is convex. 

First we reduce the computation for $FU_r$ to that for $\ov U_r \cap FU_r$.
\begin{lemma} \label{red-1}
    Assume $B = TU$ is convex. Then \[R\G_c(L\i(FU_r), \ov\BQ_\ell)[\phi] = R\G_c(L\i(\ov U_r \cap FU_r), \ov\BQ_\ell)[\phi][-2\dim (U_r \cap FU_r)]\] 
\end{lemma}
\begin{proof}
    By \cite{NTY25}, the map $(x, y) \mapsto x\i y F(x)$ gives an isomorphism \[(U_r \cap FU_r) \times (\ov U_r \cap FU_r) \overset \sim \to FU_r.\] Hence the statement follows.
\end{proof}

For $v\in\ov U_r\cap FU_r$, write $v=\prod_{\a_f\in-\Phi^+\cap F\Phi^+,\,0\le f(\bx)\le r}u_f(v_f)$ in the fixed order, and set $v_f=0$ outside this index set. For $v \in \ov U_r \cap FU_r \sm (\ov U_r \cap FU_r \cap \CK_{\phi, r})$ we set \[h(v) = \max\{f^\dag(\bx); f \in \tPhi_\re, 0 \le f(\bx) < s_f, v_f \neq 0\},\] which exists by the choice of $v$. Then there is a decomposition \[\ov U_r \cap FU_r \sm (\ov U_r \cap FU_r \cap \CK_{\phi, r}) = \bigsqcup_{h > 0} H_r^h,\] where $H_r^h = \{v \in \ov U_r \cap FU_r \sm (\ov U_r \cap FU_r \cap \CK_{\phi, r}); h(v) = h\}$.

Fix $h > 0$. Let $p: L\i(H_r^h) \to L\i(T_r H_r^h) / T_r$ be the natural projection. By Lang's theorem, $p$ is a $T_r^F$-torsor and hence \[p_! \ov\BQ_\ell = \bigoplus_{\chi \in (T_r^F)^\vee} \CE_\chi,\] where $\CE_\chi$ is a rank one local system on which $T_r^F$ acts via the character $\chi$. In particular, \[R\G_c(L^{-1}(H_r^h), \ov\BQ_\ell)[\phi]  = R\G_c(L\i(T_r H_r^h)/T_r, \CE_\phi).\] Here we put $\CE_\phi =\CE_{\phi|_{T_r^F}}$ for simplicity.

\subsection{A vanishing result} 
Let $s \in \widetilde\BR_{>0}$. Let $A_r^s$ be the set of affine roots $f \in \tPhi_\re$ with $0 \le f(\bx) \le  r$ such that $f(\bx) \ge \max\{s, s_f+\}$, or equivalently, $f(\bx) \ge \max\{s, f^\dag(\bx)+\}$. Define \[\CI_{\phi, r}^s = T_{s:r} G_r^{A_r^s}\] which is an $F$-stable subgroup of $\CI_{\phi, r}^+$. 

For $h \in \BR_{>0}$ let $B_r^h$ be the set of affine roots $f \in \tPhi_\re$ with $0 \le f(\bx) \le  r$ such that $f^\dag(\bx) \le \max\{h, s_f\}$, or equivalently, $f^\dag(\bx) \le \max\{h, f(\bx)\}$.
\begin{lemma} \label{root-closure}
    Let $f, f' \in \tPhi_\re$ with $\a_f + \a_{f'} \neq 0$.  Let $f'' = if + j f' + m_{f, f', i, j}$ be as in Lemma \ref{comm} with $f''(\bx)\le r$. Then \begin{itemize}
        \item $f'' \in B_r^h$ if $f, f' \in B_r^h$

        \item $f'' \in A_r^s$ if $f \in B_r^h$ and $f' \in A_r^s$ with $s \ge h > 0$.
    \end{itemize}
\end{lemma}
\begin{proof}
    By Lemma \ref{comm}, \begin{align*} \tag{a}
        (f'')^\dag(\bx) &=  r_{f''} - m_{f, f', i, j} - i f(\bx) - j f'(\bx) \\ &\le r_{f''} - m_{f, f', i, j} - f(\bx) - f'(\bx) \\ &\le \max\{s_{f''}, r_f, r_{f'}\} - f(\bx) - f'(\bx) \\ &\le \max\{s_{f''}-f(\bx)-f'(\bx), f^\dag(\bx) - f'(\bx),  (f')^\dag(\bx) - f(\bx)\}.
    \end{align*}

    Suppose $f, f' \in B_r^h$. Then $f^\dag(\bx) \le \max\{f(\bx), h\}$ and $(f')^\dag(\bx) \le \max\{f'(\bx), h\}$. By (a) we have \[(f'')^\dag(\bx) \le \max\{s_{f''}-f(\bx)-f'(\bx), f^\dag(\bx) - f'(\bx),  (f')^\dag(\bx) - f(\bx)\} \le \max\{s_{f''}, f''(\bx), h\},\] which means that $f'' \in B_r^h$ by noticing that $(f'')^\dag(\bx) \le s_{f''} \Leftrightarrow (f'')^\dag(\bx) \le f''(\bx)$.

    Suppose $f \in B_r^h$ and $f' \in A_r^s$ with $s \ge h$. Then $f^\dag(\bx) \le \max\{f(\bx), h\}$ and $f'(\bx) \ge \max\{s, s_{f'}+\}$. In particular, $(f')^\dag(\bx) < f'(\bx) \le f''(\bx)$. If $s_{f''} < s$, then $f''(\bx) \ge f'(\bx) \ge s = \max\{s, s_{f''}+\}$ and $f'' \in A_r^s$ as desired. Assume $s_{f''} \ge s$. It remains to show that $f''(\bx) > s_{f''}$, that is, $(f'')^\dag(\bx) < f''(\bx)$. If $f''(\bx)\le s_{f''}$, each term in the last maximum in (a) is $<s_{f''}$. By (a) we have \[(f'')^\dag(\bx) <s_{f''}\le(f'')^\dag(\bx),\] Hence a contradiction.  
\end{proof}

\begin{corollary} \label{normalize}
    Let $h \in \BR_{> 0}$ and $s \in \widetilde\BR_{\ge h}$. Then $\CI_{\phi, r}^s$ is normalized by $H_r^h$ and $H_r^h (\CI_{\phi, r}^s \cap \ov U_r \cap FU_r) = H_r^h$.
\end{corollary}
\begin{proof}
We need to show $[H_r^h,\CI_{\phi, r}^s]\subset\CI_{\phi, r}^s$.
First, let $f\in B_r^h$, $f'\in A_r^s$,  If $\a_f+\a_{f'}\ne0$, by Lemma \ref{root-closure}, we have $[G^f_r,G^{f'}_r]\subset G_r^{A_r^s}$. If $\a_{f'}=-\a_f$, we also have $[G^f_r,G^{f'}_r]\subset T_{s:r}G_r^{A_r^s}$ since $f(\bx)+f'(\bx)\ge s$ and $2f(\bx)+f'(\bx)>s_f$ and $f(\bx)+2f'(\bx)> s_f$. 
Second, for $f\in B_r^h$, every element of $[G^f_r,T_{s:r}]$ has depth $\ge f(\bx)+s$. Since $r_f-f(\bx)=f^\dag(\bx)\le\max\{f(\bx),h\}\le\max\{f(\bx),s\}$, we have $f(\bx)\ge s_f$ or $f(\bx)+s\ge r_f$. Therefore $f(\bx)+s>s_f$ and we have $[G^f_r,T_{s:r}]\subset G_r^{A_r^s}$.

By the first part of Lemma \ref{root-closure}, the subset $V_h\subseteq  V=\ov U_r\cap FU_r$ whose nonzero coordinates are indexed by $B_r^h$ is a subgroup. Since the intersection $\CI_{\phi,r}^s\cap V$ has coordinates supported in $A_r^s\subset B_r^h$, we have $\CI_{\phi,r}^s\cap V\subseteq V_h$. For $v\in H_r^h$ and $u\in\CI_{\phi,r}^s\cap V$, we therefore have $vu\in V_h$. Choose $f\in \widetilde{\Phi}_{\mathrm{re}}$ such that $v_f\ne0,0\le f(\bx)<s_f,f^\dag(\bx)=h$. Then $f(\bx)=r_f-h<h$. So $u\in G_{s:r}\subset G_{h:r}$ implies that $(vu)_f\ne0$. Therefore $h(vu)=h$. We see that $H_r^h(\CI_{\phi,r}^s\cap V)\subset H_r^h$. By taking $u=1$ we get equality.
\end{proof}

We fix $g \in G_r$ with $L(g) \in H_r^h$. Consider the morphism \[\pi: G_r \to G_r, \quad x \mapsto x\i L(g)F(x).\] Let $q: L\i(T_r H_r^h)/ T_r \to  G_r / (T_r\CI_{\phi, r}^h)$ be the natural quotient map. By Corollary \ref{normalize}, the left multiplication by $g\i$ gives an identification \[ q\i(g T_r \CI_{\phi, r}^h) \cong \{x \in G_r^{A_r^h}; \pi(x) \in T_r H_r^h\}.\] Using the above isomorphism, we define \[\pi_{T_{r}}: q\i(g T_r\CI_{\phi, r}^h) \to T_r, \quad x \mapsto \pi(x)_{T_r}.\] Here $x_{T_r}$ denotes the torus part of $x \in F(T_r U_r)$.

Let $\D_r^h = \{f \in A_r^h; \a_f \in \Phi^+ \cap -F\Phi^+\}$ and $D_r^h = \{f \in \D_r^h; f(\bx) = h\}$.
\begin{proposition} \label{fiber}
    Assume $B = TU$ is convex. Let $h > 0$ and $g \in L\i(H_r^h)$. Write $y=L(g)\i$. Then the map $x \mapsto (x_f)_{f \in \D_r^h}$ induces an isomorphism \[q\i(gT_r\CI_{\phi, r}^h) \cong \BA^{\D_r^h}.\] Moreover, under the above identification, \[\pi_{T_r}(x) \equiv \prod_{f \le f' \in D_r^h} \l_{\a_f^\vee, r_f} \circ \mu_{c_{f, f'} y_{f^\dag}} (x_{f'}) \mod \hCI_{\phi, r}^\dag,\] where $c_{f, f'} \in \ov\BF_q$ are constants depending on $y$ with $c_{f, f} =  1$ for $f \in D_r^h$.
\end{proposition}
\begin{proof}
    It follows by applying \cite[Proposition {3.6}]{IN25} and the proof of \cite[Proposition {5.6}]{IN25}.
\end{proof}

\begin{corollary} \label{vanish-CE}
    We have $R\G_c(q\i(g T_r\CI_{\phi, r}^h), \CE_\phi) = 0$.
\end{corollary}
\begin{proof}
    By \cite[Lemma 5.2]{IN25}, we have \[\CE_\phi|_{q\i(g T_r\CI_{\phi, r}^h)} \cong \pi_{T_r}^* \CL_\phi.\] By Proposition \ref{fiber} and Lemma \ref{L_phi}, we have \[\pi_{T_r}^* \CL_\phi \cong \eta^*\CL_\phi,\] where the second isomorphism follows from Lemma \ref{L_phi} and \[\eta(x) =\prod_{f \le f' \in D_r^h} \l_{\a_f^\vee, r_f} \circ \mu_{c_{f, f'} y_{f^\dag}} (x_{f'}),\] with $c_{f, f} = 1$ for $f \in D_r^h$. Let $f_0 = \min\{f \in D_r^h; y_{f^\dag} \neq 0\}$ which exists by the choice of $g$. Then \[\eta(x) = \l_{\a_{f_0}^\vee, r_{f_0}} \circ \mu_{y_{f_0^\dag}}(x_{f_0})  \prod_{f \le f' \in D_r^h,\ f' > f_0} \l_{\a_f^\vee, r_f} \circ \mu_{c_{f, f'} y_{f^\dag}}(x_{f'}).\] By Lemma \ref{L_phi} and \cite[Proposition 5.7]{IN25}, we have $R\G_c(q\i(gT_r\CI_{\phi, r}^h), \eta^*\CL_\phi) = 0$ as desired.   
\end{proof}

\begin{corollary}
     Assume $B=TU$ is convex. Then \[R\G_c(L\i(\ov U_r \cap FU_r \cap \CK_{\phi, r}), \ov\BQ_\ell)[\phi] = R\G_c(L\i(\ov U_r \cap FU_r), \ov\BQ_\ell)[\phi].\] 
\end{corollary}
\begin{proof}
     By Corollary \ref{vanish-CE} we have \[R\G_c(L\i(H_r^h), \ov\BQ_\ell)[\phi] =  R\G_c(L\i(T_r H_r^h)/T_r, \CE_\phi) = 0 \text{ for } h > 0.\] Hence the statement follows.
\end{proof}

\subsection{The reduction to $F\CI_{\phi, r}^+$}
Note that $\ov U_r \cap FU_r \cap \CK_{\phi, r} \subseteq F \CI_{\phi, r}^+$.
\begin{lemma} \label{cartesian}
     There is a Cartesian square
     \[ 
     \xymatrix@C=5em{
     L\i(\ov U_r \cap FU_r \cap \CK_{\phi, r}) \times T_{0+:r} \ar[d]_\mu  \ar[r]^{\quad\quad\quad\quad\quad \pr_2} & T_{0+:r} \ar[d]^{F\i L} \\
     L\i(T_{0+:r}(\ov U_r \cap FU_r \cap \CK_{\phi, r})) \ar[r]^{\quad\quad\quad\quad\quad \t_r \circ F\i \circ L} & T_{0+:r},}
     \] where $\mu$ is given by $(x, y) \mapsto xy$. In particular, we have \begin{align*} &\quad\ R\G_c( L\i(\ov U_r \cap FU_r \cap \CK_{\phi, r}), \ov\BQ_\ell)[\phi] \\ &\cong R\G_c( L\i(T_{0+:r}(\ov U_r \cap FU_r \cap \CK_{\phi, r})), (F\i L)^*\CF)[\phi_0][2 M_0],\end{align*} where $M_0 = \dim T_{0+:r}$.
\end{lemma}
\begin{proof}
    The Cartesian square follows by definition. By the base change theorem we have \begin{align*}
        &\quad\ R\G_c( L\i(\ov U_r \cap FU_r \cap \CK_{\phi, r}), \ov\BQ_\ell)[\phi|_{T_{0+:r}^F}] \\ &\cong R\G_c( L\i(\ov U_r \cap FU_r \cap \CK_{\phi, r}) \times T_{0+:r}, \ov\BQ_\ell)[\phi|_{T_{0+:r}^F}] [2M_0] \\ &\cong R\G_c( L\i(T_{0+:r} (\ov U_r \cap FU_r \cap \CK_{\phi, r})), (\t_r F\i L)^*\CL) [2M_0]\\ &\cong R\G_c( L\i(T_{0+:r} (\ov U_r \cap FU_r \cap \CK_{\phi, r})), (F\i L)^*\CF) [2M_0].
    \end{align*} So the statement follows.
\end{proof}

\begin{proposition} \label{new-constr}
    Assume $B= TU$ is convex. Then \[R\G_c(L\i(\ov U_r \cap FU_r \cap \CK_{\phi, r}), \ov\BQ_\ell)[\phi] = R\G_c(L\i(F\CI_{\phi, r}^+), (F\i L)^*\CF)[\phi_0][2M],\] where $M = \dim T_{0+:r} + \dim (U_r \cap FU_r \cap \CK_{\phi, r}) + \dim (\ov U_r \cap \CK_{\phi, r}^+)$.

    In particular, $H_c^*(L\i(\ov U_r \cap FU_r \cap \CK_{\phi, r}), \ov\BQ_\ell)[\phi]$ is independent of the choice of $r \ge r_\phi$.
\end{proposition}
\begin{proof}
    As $B$ is convex, the statement follows from \cite[Proposition 6.3]{IN25} that the map $(x, y, z) \mapsto xyz$ induces an isomorphism \[L\i(T_{0+:r}(\ov U_r \cap FU_r \cap \CK_{\phi, r})) \times (\ov U_r \cap \CK_{\phi, r}^+) \times (U_r \cap FU_r \cap \CK_{\phi, r})  \overset \sim \to L\i(F\CI_{\phi, r}^+).\] By Lemma \ref{local-CF}, Lemma \ref{cartesian} and the projection formula we have \begin{align*}
        &\quad\ R\G_c( L\i(\ov U_r \cap FU_r \cap \CK_{\phi, r}),\ov\BQ_\ell)[\phi] \\ &\cong  R\G_c( L\i(T_{0+:r}(\ov U_r \cap FU_r \cap \CK_{\phi, r})), (F\i L)^*\CF)[\phi_0][2 M_0] \\ &\cong R\G_c( L\i(F \CI_{\phi, r}^+), (F\i L)^*\CF)[\phi_0][2 M].
    \end{align*} Hence the statement follows.
\end{proof}

\section{The intertwining space} \label{sec:inter}

Let $w \in G^F$. For a virtual smooth $G_{\bx, 0}^F$-module $\rho$, the virtual intertwining space is defined by \[I_w(\rho) = \Hom_{G_{\bx, 0}^F \cap {}^w G_{\bx, 0}^F}(\rho, {}^w\rho) = (\rho^\vee \otimes {}^w \rho)^{G_{\bx, 0}^F \cap {}^w G_{\bx, 0}^F}.\]

Let $r'\gg r \ge r_\phi$ be a sufficiently large integer such that ${}^{w\i} G_{\bx, r'+} \subseteq G_{\bx, r+}$. Consider the variety \[\Sigma = \{(x_1, x_2, y) \in \CI_{\phi, r'}^+ \times \CI_{\phi, r}^+ \times G_{\bx, 0} w G_{\bx, 0} / G_{\bx, r+}; F(x_1 y x_2\i) = y\},\] on which $T_0^F \times T_0^F$ acts by $(t_1, t_2): (t_1 x_1 t_1\i, t_2 x_2 t_2\i, t_1 y t_2\i)$.

Let $\CM = \pr_{12}^*(\CF \boxtimes \CF^\vee)$. Here $\pr_{12}: \Sigma \to \CI_{\phi, r'}^+ \times \CI_{\phi, r}^+$ is the natural projection given by $(x_1, x_2, y) \mapsto (x_1, x_2)$.
\begin{proposition} \label{inter-CM}
    We have \[ \dim I_w(R_{T_r}^{G_r}(\phi)) = \dim H_c^*(\Sigma, \CM)[\phi_0 \boxtimes \phi_0\i].\]
\end{proposition}
\begin{proof}
   Let $P_{r', r} \subseteq G_{r'} \times G_r$ be the graph of the homomorphism $(G_{\bx,0}\cap{}^wG_{\bx,0})/G_{\bx,r'+} \to G_r$ given by $g_1 \mapsto w\i g_1 w$. By \cite[Proposition 4.5]{LNY26}, there is a $P_{r', r}^F$-torsor \[\th: L\i(F \CI_{\phi, r'}^+) \times L\i(F \CI_{\phi, r}^+) \to \Sigma, \quad (g_1, g_2) \mapsto (F\i L(g_1), F\i L(g_2), g_1\i w g_2).\] By definition $\th^* \CM = (F\i L)^*\CF \boxtimes (F\i L)^*\CF^\vee$. By the projection formula we have \[\th_! \th^* \CM \cong \th_! \ov\BQ_\ell \otimes \CM.\] As $\th$ is a $P_{r', r}^F$-torsor, $(\th_! \ov\BQ_\ell \otimes \CM)^{P_{r' ,r}^F} \cong (\th_! \ov\BQ_\ell)^{P_{r' ,r}^F} \otimes \CM \cong \CM$. Hence \[R\G_c(\Sigma, \CM) \cong R\G_c(L\i(F \CI_{\phi, r'}^+) \times L\i(F \CI_{\phi, r}^+), (F\i L)^*\CF \boxtimes (F\i L)^*\CF^\vee)^{P_{r', r}^F}.\] As the action of $T_0^F \times T_0^F$ on $L\i(F \CI_{\phi, r'}^+) \times L\i(F \CI_{\phi, r}^+)$ by right multiplication commutes with the action of $P_{r', r}^F$ by left multiplication, it follows that  \begin{align*}&\quad\ H_c^*(\Sigma, \CM)[\phi_0 \boxtimes \phi_0\i] \\ &\cong H_c^*(L\i(F \CI_{\phi, r'}^+) \times L\i(F \CI_{\phi, r}^+), (F\i L)^*\CF \boxtimes (F\i L)^*\CF^\vee)[\phi_0 \boxtimes \phi_0\i]^{P_{r', r}^F} \\  &\cong  (H_c^*( L\i(F \CI_{\phi, r'}^+), (F\i L)^*\CF)[\phi_0] \otimes  H_c^*( L\i(F \CI_{\phi, r}^+), (F\i L)^*\CF^\vee)[\phi_0\i])^{P_{r', r}^F} \\ &\cong (R_{T_r}^{G_r}(\phi) \otimes R_{T_r}^{G_r}(\phi)^\vee)^{P_{r', r}^F} \\ &\cong I_w(R_{T_r}^{G_r}(\phi))^\vee,\end{align*} where the third isomorphism follows from Proposition \ref{new-constr}. The proof is finished.
\end{proof}

\section{Vanishing results on $\CM$} \label{sec:Sigma}
Let the notation be as in previous sections. Let $N_T$ be the normalizer of $T$ in $G$. Denote by \[W = N_T(\brk) / T(\brk) \text{ and } \tW = N_T(\brk) / T(\CO_\brk) = X_*(T) \rtimes W\] the Weyl group and the Iwahori-Weyl group of $G$ respectively. Recall the Bruhat decompositions \[G = \bigsqcup_{x \in \tW} \CI x \CI, \quad  G_{\bx, 0} = \bigsqcup_{x \in W_\bx} \CI x \CI,\] where $W_\bx \subseteq \tW$ denotes the Weyl group of $G_{\bx, 0}$, and $\CI$ is the Iwahori subgroup generated by $T(\COk)$ and $G^f$ for $f \in \tPhi_\re^{>0}$.

For $w \in G^F$ we have \[G_{\bx, 0} w G_{\bx, 0} / G_{\bx, r+} = \bigsqcup_{z \in W_\bx z_w W_\bx} C_z,\] where $C_z =  \CI z \CI / G_{\bx, r+}$ and $z_w \in \tW$ satisfies $w \in \CI z_w \CI$. This induces a decomposition \[\Sigma = \bigsqcup_{z \in W_\bx z_w W_\bx} \Sigma_z,\] where $\Sigma_z = \{(x_1, x_2, y) \in \CI_{\phi, r'}^+ \times \CI_{\phi, r}^+ \times C_z; F(x_1 y x_2\i) = y\}$.

\subsection{A deformation method}
We fix $z \in W_\bx z_w W_\bx$. Define \begin{align*}
    A &=\{f \in \tPhi_\re^{>0}; \text{ $z\i(f)(\bx) \le r$, either $z\i(f)(\bx) < f(\bx)$ or $z\i(f)(\bx) = f(\bx)$ and $\a_f > 0$}\}; \\ B &=\{f \in \tPhi_{\re, r}^{>0}; \text{ either $z(f)(\bx) < f(\bx)$ or $z(f)(\bx) = f(\bx)$ and $\a_{z(f)} < 0$}\}.
\end{align*} We have the following parameterization result.
\begin{lemma}[{\cite[Lemma 5.2]{LNY26}}] \label{coordinate}
    There is an isomorphism \[ G_{r'}^A \times zT_r \times G_r^B \to C_z, \quad\ (v_1, \dz, v_2) \mapsto v_1 \dz v_2\i.\]
\end{lemma}

Let $C_z^0$ be the set of points $y = v_1 \dz v_2\i \in C_z$ as in Lemma \ref{coordinate} such that $v_1 \in \CI_{\phi, r'}^+$ and $v_2 \in \CI_{\phi, r}^+$, that is, $(v_1)_{f_1} = 0 = (v_2)_{f_2}$ for all $f_1 \in A$ and $f_2 \in B$ such that $f_1 < s_{f_1}$ and $f_2 < s_{f_2}$. Let $\Sigma_z^0$ be the subset of triples $(x_1, x_2, y) \in \Sigma_z$ such that $y \in C_z^0$.

For $y \in C_z$ we denote by \[\Sigma_z(y) = \{(x_1, x_2) \in \CI_{\phi, r'}^+ \times \CI_{\phi, r}^+; F(x_1 y x_2\i) = y\}\] the fiber of the projection map $\Sigma_z \to C_z$ given by $(x_1, x_2, y) \mapsto y$.

We fix $y = v_1 \dz v_2\i \in C_z \sm C_z^0$ with $v_1 \in G_{r'}^A$ and $v_2 \in G_r^B$. Let \[A'= \{f \in A; (v_1)_f \neq 0, f < s_f\}, \quad B'= \{f \in B; (v_2)_f \neq 0, f < s_f\}.\] Define $f_1 \in A'$ (resp. $f_2 \in B'$) such that $f_1^\dag = \max\{f^\dag; f \in A'\}$ if $A' \neq \emptyset$ (resp. $f_2^\dag = \max\{f^\dag; f \in B'\}$ if $B' \neq \emptyset$) and $f_1 = 0$ (resp. $f_2 = 0$) otherwise.
By the choice of $y$, we have $f_1 \neq 0$ or $f_2 \neq 0$.

\begin{lemma} \label{deform}
    Assume $z\i(f_1^\dag) > s_{z\i(f_1)}, f_2^\dag$ (resp. $z(f_2^\dag) > s_{z(f_2)}, f_1^\dag$) or $r_{z(f_2)} > r_{f_2}, f_1^\dag$ (resp. $r_{z\i(f_1)} > r_{f_1}, f_2^\dag$ ). Then there exists an isomorphism \[\D: \BA^1 \times \Sigma_z(y) \to \BA^1 \times \Sigma_z(y), \quad (\t, x_1, x_2) \mapsto (\t, x_1 \d_1(\t), x_2 \d_2(\t))\] such that \[\D^*(\ov\BQ_\ell \boxtimes \CM) \cong \CP\boxtimes\CM,\] where $\CP = \d_1^*\CF \otimes \d_2^*\CF^\vee$ is a nontrivial Artin-Schreier sheaf. In particular, $R\G_c(\Sigma_z(y), \CM) = 0$.
\end{lemma}
\begin{proof}
    We define \[\d_1 = v_1 \d v_1\i, \quad \d_2 = v_2 {}^{\dz{\i}} \d v_2\i,\] where $\d: \BA^1 \to \CI_{\phi, r'}^+$ is a suitable morphism such that \[\tag{$\ast$} \d_1:\BA^1\to\CI_{\phi,r'}^+,\quad\d_2:\BA^1\to\CI_{\phi,r}^+, \quad \CP = \d_1^*\CF \otimes \d_2^*\CF^\vee \cong \l_{\a^\vee, r_\a}^*\CF\] for some root $\a \in \Phi$ with $r_\a>0$. The map $(\t, x_1, x_2) \mapsto (\t, x_1 \d_1, x_2 \d_2)$ is an automorphism $\D$ of $\BA^1 \times \Sigma_z(y)$. Moreover, by Lemma \ref{L_phi} and Lemma \ref{local-CF} we have $\D^*(\ov\BQ_\ell \boxtimes \CM) \cong \CP\boxtimes\CM$ with $\CP$ a non-trivial Artin-Schreier local system. Hence \begin{align*} &\quad\ R\G_c(\Sigma_z(y), \CM)[-2] \\ &\cong R\G_c(\BA^1 \times \Sigma_z(y), \ov\BQ_\ell \boxtimes \CM) \\ &\cong R\G_c(\BA^1 \times \Sigma_z(y), \D^*(\ov\BQ_\ell \boxtimes \CM)) \\ &\cong R\G_c(\BA^1 \times \Sigma_z(y), \CP \boxtimes \CM) \\ &\cong R\G_c(\BA^1, \CP) \otimes R\G_c(\Sigma_z(y), \CM) \\ &=0. \end{align*} Therefore, it remains to construct $\d$ such that $(*)$ is satisfied.

    Case (1): $z\i(f_1^\dag) > s_{z\i(f_1)}, f_2^\dag$. Take $\d(\t) = u_{f_1^\dag}((v_1)_{f_1}\i \t)$. Then $\d \equiv1\;\mod \hCI_{\phi, r'}^\dag$ and ${}^{\dz\i} \d \equiv1\mod \hCI_{\phi, r}^\dag$ since $f_1^\dag > f_1$ and $z\i(f_1^\dag) > s_{z\i(f_1)}$ respectively. As $z\i(f_1^\dag) > f_2^\dag$, it follows from Lemma \ref{commutator} that \[\d_2 = v_2 {}^{\dz\i} \d v_2\i = [v_2, {}^{\dz\i} \d] {}^{\dz\i} \d   \equiv1\mod \hCI_{\phi, r}^\dag\] and \[\d_1 = v_1 \d v_1\i = \d [\d\i, v_1] \equiv \l_{\a_{f_1}^\vee, r_{f_1}} \mod \hCI_{\phi, r'}^\dag.\] Hence $(*)$ holds as desired.

    Case (2): $z(f_2^\dag) > s_{z(f_2)}, f_1^\dag$. Take $\d(\t) = {}^\dz u_{f_2^\dag}((v_2)_{f_2}\i \t)$. Then ${}^{\dz\i}\d \equiv1\mod \hCI_{\phi, r}^\dag$ and $\d \equiv1\mod \hCI_{\phi, r'}^\dag$ since $f_2^\dag > f_2$ and $z(f_2^\dag) > s_{z(f_2)}$ respectively. As $z(f_2^\dag) > f_1^\dag$, it follows from Lemma \ref{commutator} that \[\d_1 = v_1 \d v_1\i = [v_1, \d] \d\equiv1\mod \hCI_{\phi, r'}^\dag\] and \[\d_2 = v_2 {}^{\dz\i}\d v_2\i = {}^{\dz\i}\d [{}^{\dz\i}\d\i, v_2] \equiv \l_{\a_{f_2}^\vee, r_{f_2}} \mod \hCI_{\phi,r}^\dag.\] Hence $(*)$ holds as desired.

    Case (3): $r_{z(f_2)} > r_{f_2}, f_1^\dag$ (resp.  $r_{z\i(f_1)} > r_{f_1}, f_2^\dag$). Take $\d = \l_{\a_{z(f_2)}^\vee, r_{z(f_2)}}$ (resp. $\d = {}^\dz \l_{\a_{z\i(f_1)}^\vee, r_{z\i(f_1)}}$). As $r_{z(f_2)} > r_{f_2}$ (resp. $r_{z\i(f_1)} > r_{f_1}$), we have ${}^{\dz\i}\d = \l_{\a_{f_2}^\vee, r_{z(f_2)}} \equiv1\mod \hCI_{\phi, r}^\dag$ (resp. $\d = \l_{\a_{f_1}^\vee, r_{z\i(f_1)}} \equiv1\mod \hCI_{\phi, r'}^\dag$). As $r_{z(f_2)} > f_1^\dag$ (resp. $r_{z\i(f_1)} > f_2^\dag$), it follows by Lemma \ref{commutator} that \begin{align*} &\d_1 = v_1 \d v_1\i \equiv \l_{\a_{z(f_2)}^\vee, r_{z(f_2)}} \mod \hCI_{\phi, r'}^{\dag} \\ (\text{resp. } & \d_2 = v_2 {}^{\dz\i}\d v_2\i \equiv \l_{\a_{z\i(f_1)}^{\vee}, r_{z\i(f_1)}} \mod \hCI_{\phi, r}^\dag). \end{align*} Moreover, as $r_{z(f_2)} > r_{f_2} > f_2^\dag$ (resp. $r_{z\i(f_1)} > r_{f_1} > f_1^\dag$), by Lemma \ref{commutator} we have \[\text{$\d_2 = v_2 {}^{\dz\i}\d v_2\i  =  [v_2, {}^{\dz\i}\d] {}^{\dz\i}\d \equiv1\mod \hCI_{\phi, r}^\dag$ (resp. $\d_1 = v_1 \d v_1\i = [v_1, \d] \d \equiv1\mod\hCI_{\phi, r'}^\dag$)}.\] Hence $(*)$ holds as desired.
\end{proof}

\subsection{Reduction to the locus $\Sigma_z^0$}
We recall the following two lemmas. 
\begin{lemma} [{\cite[Lemma 5.1]{LNY26}}] \label{value}
    Let $f \in \tPhi_\re$. We have $z(f)(\bx) - f(\bx) = \a_f(z\i(\bx) - \bx)$. In particular, $z(f)(\bx) \le f(\bx) \Leftrightarrow z(f^\dag)(\bx) \ge f^\dag(\bx)$.
\end{lemma}

\begin{lemma} [{\cite[Lemma 5.4]{LNY26}}] \label{gedag}
    If $f_1\ne0, r_{f_1} \ge r_{z\i(f_1)}$ (resp. $f_2\ne0, r_{f_2} \ge r_{z(f_2)}$), then $z\i(f_1^\dag) > s_{z\i(f_1)}$ (resp. $z(f_2^\dag) > s_{z(f_2)}$). 
\end{lemma}  

The main result of this section is as follows.
\begin{proposition} \label{vanishing}
    Let $y \in C_z$. Then $R\G_c(\Sigma_z(y), \CM) = 0$ unless $y \in C_z^0$. In particular, we have $R\G_c(\Sigma_z^0, \CM) \cong R\G_c(\Sigma_z, \CM)$.
\end{proposition}
\begin{proof}
    It follows in a similar way as \cite[Proposition 5.3]{LNY26}. We include a proof for completeness.

    By Lemma \ref{value} and the definitions of $A$ and $B$, we have $z(f_2^\dag)(\bx) \ge f_2^\dag(\bx)$ and $z\i(f_1^\dag)(\bx) \ge f_1^\dag(\bx)$.

    Case (1): either $f_1^\dag(\bx) < f_2^\dag(\bx)$ or $f_1^\dag(\bx) = f_2^\dag(\bx)$ and $z(f_2^\dag)(\bx) > f_2^\dag(\bx)$. Hence $z(f_2^\dag)(\bx) > f_1^\dag(\bx)$ and $z(f_2^\dag) > f_1^\dag$. If $r_{z(f_2)} > r_{f_2}$, then $r_{z(f_2)} > f_1^\dag$ (since $r_{z(f_2)} > r_{f_2} \ge f_2^\dag(\bx) \ge f_1^\dag(\bx)$), and the statement follows by Lemma \ref{deform}. So we can assume that $r_{f_2} \ge r_{z(f_2)}$. By Lemma \ref{gedag} we have $z(f_2^\dag) > s_{z(f_2)}$. The statement also follows by Lemma \ref{deform}.

    Case (2): either $f_1^\dag(\bx) > f_2^\dag(\bx)$ or $f_1^\dag(\bx) = f_2^\dag(\bx)$ and $z\i(f_1^\dag)(\bx) > f_1^\dag(\bx)$. Hence $z\i(f_1^\dag)(\bx) > f_2^\dag(\bx)$ and $z\i(f_1^\dag) > f_2^\dag$. If $r_{z\i(f_1)} > r_{f_1}$, then $r_{z\i(f_1)} > f_2^\dag$ (since $r_{z\i(f_1)} > r_{f_1} \ge f_1^\dag(\bx) \ge f_2^\dag(\bx)$), and the statement follows by Lemma \ref{deform}. So we can assume that $r_{f_1} \ge r_{z\i(f_1)}$. By Lemma \ref{gedag} we have $z\i(f_1^\dag) > s_{z\i(f_1)}$. The statement also follows by Lemma \ref{deform}.

    Case (3): $z\i(f_1^\dag)(\bx) = f_1^\dag(\bx) = f_2^\dag(\bx) = z(f_2^\dag)(\bx)$. By Lemma \ref{value}, we have $z\i(f_1)(\bx) = f_1(\bx)$, which means  $\a_{f_1} > 0$ by the definition of $A$. Similarly we have $z(\a_{f_2}) < 0$. Therefore, $\a_{z(f_2^\dag)} > 0 > \a_{f_1^\dag}$ and hence $z(f_2^\dag) > f_1^\dag$. If $r_{f_2} < r_{z(f_2)}$, then $r_{z(f_2)} > f_1^\dag$ (since $r_{z(f_2)} > r_{f_2} \ge f_2^\dag(\bx) = f_1^\dag(\bx)$) and the statement follows from Lemma \ref{deform}. So we can assume that $r_{f_2} \ge r_{z(f_2)}$. By Lemma \ref{gedag} we have $z(f_2^\dag) > s_{z(f_2)}$. The statement follows again by Lemma \ref{deform}.
\end{proof}

\section{The dimension formula} \label{sec:dim-inter}
In this section, we deduce the dimension formula of $I_w(R_{T_r}^{G_r}(\phi))$. By Lemma \ref{coordinate} we have a natural isomorphism \[G_{r'}^{A \cap E_\phi} \times z T_r \times G_r^{B \cap E_\phi} \cong C_z^0, \quad (v_1, \dz, v_2) \mapsto v_1 \dz v_2\i.\] Using the above identification, we consider the variety \[\Sigma_z^{0, \flat} = \{(x_1, x_2, v_1 \dz v_2\i) \in \CI_{\phi, r'}^+ \times \CI_{\phi, r}^+ \times C_z^0; F(x_1 \dz x_2\i) = v_1 \dz v_2\i\},\] on which $T_{0}^F \times T_{0}^F$ acts by \[(t_1, t_2): (x_1, x_2, v_1, v_2, \dz) \mapsto ({}^{t_1} x_1, {}^{t_2} x_2, {}^{t_1} v_1, {}^{t_2} v_2, t_1 \dz t_2\i).\] 

We denote by $\CM^\flat$ the pull-back of $\CF \boxtimes \CF^\vee$ under the natural projection $\Sigma_z^{0, \flat} \to \CI_{\phi, r'}^+ \times \CI_{\phi, r}^+$.
\begin{lemma} \label{flat}
    The isomorphism \[f: \Sigma_z^0 \to \Sigma_z^{0, \flat}, \quad (x_1, x_2, v_1 \dz v_2\i) \mapsto (x_1 v_1, x_2 v_2, v_1 \dz v_2\i)\] induces a $T_0^F \times T_0^F$-equivariant isomorphism of sheaves \[f^*\CM^\flat \cong \CM.\] In particular, $R\G_c(\Sigma_z^0, \CM)[\phi_0 \boxtimes \phi_0\i] \cong R\G_c(\Sigma_z^{0, \flat}, \CM^\flat)[\phi_0 \boxtimes \phi_0\i]$.
\end{lemma}
\begin{proof}
    It follows from Lemma {\color{red}\ref{local-CF}} by noticing that $\t_{r'}(v_1)=\t_r(v_2)=1$.
\end{proof}

Consider a cover of $\Sigma_z^{0, \flat}$ defined by \[\hSigma_z^{0, \flat} = \{((x_1, x_2, y), s_1, s_2) \in \Sigma_z^{0, \flat} \times T_{0+:r'} \times T_{0+:r}; (\t_{r'}(x_1), \t_r(x_2)) = (L(s_1), L(s_2))\},\] on which $T_{r'}^F \times T_r^F$ acts by \[(t_1 \t_1, t_2 \t_2): (\xi, s_1,s_2) \mapsto ((t_1, t_2) \cdot \xi, s_1\t_1\i, s_2\t_2\i),\] where $\xi \in \Sigma_z^{0, \flat}$, $(t_1, t_2) \in T_0^F \times T_0^F$ and $(\t_1, \t_2) \in T_{0+:r'}^F \times T_{0+:r}^F$. 
\begin{lemma} \label{cover}
    We have $R\G_c(\Sigma_z^{0, \flat}, \CM^\flat)[\phi_0 \boxtimes \phi_0\i] \cong R\G_c(\widehat\Sigma_z^{0, \flat}, \ov\BQ_\ell)[\phi \boxtimes \phi\i]$.
\end{lemma}
\begin{proof}
By construction there is a Cartesian diagram
\[ 
     \xymatrix{
     \hSigma_z^{0, \flat}  \ar[d]_{\pr_1}  \ar[r]^{\pr_{23} \quad\quad\ } & T_{0+:r'} \times T_{0+:r} \ar[d]^{L \times L} \\
     \Sigma_z^{0, \flat} \ar[r]^{f \quad\quad\ } & T_{0+:r'} \times T_{0+:r},}
     \] where $f$ is given by $(x_1, x_2, y) \mapsto (\t_{r'}(x_1), \t_r(x_2))$. Then the statement follows from the base change theorem and the identity $\CM^\flat \cong f^*(\CL \boxtimes \CL^\vee)$.
\end{proof}

\begin{proposition} \label{chi-flat}
    We have $H^*_c(\hSigma_z^{0, \flat}, \ov\BQ_\ell)[\phi \boxtimes \phi\i] \cong \ov\BQ_\ell$ if $F(z) = z$ and ${}^z \phi|_{T_r^F} = \phi |_{T_r^F}$ and $H^*_c(\hSigma_z^{0,\flat}, \ov\BQ_\ell)[\phi \boxtimes \phi\i] = 0$ otherwise.
\end{proposition}
\begin{proof}
    Let $H = \{(t_1, t_2) \in T_0 \times T_0; L(t_1) \dz L(t_2)\i = \dz\}$. Using Teichm\"{u}ller lift, we view $H$ and its identity component $H^\circ$ as subgroups of $T_{r'} \times T_r$. Note that $H$ acts on $\hSigma_z^{0, \flat}$ by \[(t_1, t_2): (x_1, x_2, v_1, v_2, \dz,s_1,s_2) \mapsto({}^{t_1} x_1, {}^{t_2} x_2, {}^{F(t_1)} v_1, {}^{F(t_2)} v_2, F(t_1) \dz F(t_2)\i,s_1,s_2).\] The actions by $H$ and by $T_{r'}^F \times T_r^F$ commute with each other. Hence \[H^*_c(\hSigma_z^{0, \flat}, \ov\BQ_\ell)[\phi \boxtimes \phi\i] \cong H^*_c((\hSigma_z^{0, \flat})^{H^\circ}, \ov\BQ_\ell)[\phi \boxtimes \phi\i].\] As the two projections $H^\circ \to T_0$ to the first and the second factors are both surjective, it follows that $(\hSigma_z^{0, \flat})^{H^\circ} = \emptyset$ unless $F(z) = z \in \tW$. So we may assume $z = F(z) \in \tW$ and choose an $F$-fixed representative $\dot z_0\in N_T(k)$. Then \begin{align*}(\hSigma_z^{0, \flat})^{H^\circ} &\cong \{(s_1, s_2, t) \in T_{0+:r'} \times T_{0+:r} \times T_r; F(L(s_1)\dot z_0tL(s_2)\i)=\dot z_0t\} \\ &\cong (\dot z_0 T_r^F) \times T_{0+:r'} \times T_{0+:r}.\end{align*} The second map sends $(s_1,s_2,t)$ to $(F(s_1)\dot z_0tF(s_2)\i,s_1,s_2)$. In particular, \[H^*_c(\hSigma_z^{0, \flat}, \ov\BQ_\ell)[\phi \boxtimes \phi\i] \cong H^*_c(\dot z_0 T_r^F, \ov\BQ_\ell)[\psi \boxtimes \psi\i].\] Here the character $\psi$ is given by $\psi(t_0\tau)=\phi(t_0\tau\i)$ for $t_0\in T_0^F,\tau\in T_{0+:r}^F$. Hence the statement follows. 
\end{proof}

\begin{theorem}
    Assume $B = T U$ is convex. For $w \in G^F$ we have \[\dim I_w(R_{T_r}^{G_r}(\phi)) = \sharp\{z \in (W_\bx z_w W_\bx)^F; {}^z \phi|_{T_r^F} = \phi |_{T_r^F}\}.\]
\end{theorem}
\begin{proof}
    By Proposition \ref{inter-CM}, Proposition \ref{vanishing}, Lemma \ref{flat}, Lemma \ref{cover} and Proposition \ref{chi-flat}.
\end{proof}

\section{Vanishing of $V_r^F$-invariants} \label{sec:sum}
As $G$ splits over $\brk$, there exists a maximally $k$-split maximal torus $S \subseteq G$ such that $\bx \in \CA(S, \brk)$. Let $\l$ be a $k$-rational non-central cocharacter of $S$. Let \[P = P_G(\l) := \{g: \lim_{t \to 0} \l(t) g \l(t)\i \text{ exists }\}\] be the associated proper parabolic subgroup of $G$. Write $P = M V$ where $M$ is the Levi subgroup containing $S$ and $V$ is the unipotent radical of $P$. In particular, there exists $\l = F(\l) \in X_*(S)$ such that the sets of roots for $S$ in $P$ and $M$ are given by $\Phi(P, S) = \{\a \in \Phi(G, S); \a(\l) \ge 0\}$ and $\Phi(M, S) = \{\a \in \Phi(G, S); \a(\l) = 0\}$ respectively. Let $\ov V$ be the unipotent radical of the opposite of $P$. We denote by $\tPhi_{\ov V}$ the set of affine roots for $S$ whose finite parts appear in $\ov V$. 

Let $W_{M, r} = N_M(S)_r / S_r$ be the Weyl group of $M_r$, where $N_{M}(S)$ is the normalizer of $S$ in $M$. Let $N_G(S, T) = \{x \in G; {}^x S = T\}$ and $W_{S, T, r} = T_r \backslash N_G(S, T)_r / S_r$.  Then we have a decomposition \[G_r = \bigsqcup_{z \in W_{S, T, r} / W_{M, r}} \CI_{z, r} z P_r,\] where $\CI_{z, r} = \CI_r \cap {}^z \ov V_r$. 

Let $z \in W_{S, T, r} / W_{M, r}$. We set $\tPhi_{z, r} = \{f \in \tPhi_{\re, r}^{>0}; z\i(f) \in \tPhi_{\ov V}\}$. Then \[G_r^{\tPhi_{z, r}} = \CI_{z, r}.\] For $v \in \CI_{z, r}$ we denote by $v_f \in \BA_f = \ov\BF_q$ the coordinate of $v$ at $f \in \tPhi_{z, r}$ via the above identification.

\subsection{Reduction to the locus $\Xi_z^0$}
Consider the variety \[\Xi = \{(x, \psi) \in \CI_{\phi, r}^+ \times G_r / V_r; F(x\psi) = \psi\},\] on which $T_0^F$ acts by $t: (x, \psi) \mapsto (t\i x t, t\i \psi)$.
\begin{lemma} \label{V-geo}
    Assume $B = T U$ is convex. Then $R_{T_r}^{G_r}(\phi)^{V_r^F} \cong H_c^*(\Xi, \pr_1^*\CF)[\phi_0]$.
\end{lemma}
\begin{proof}
    It follows in a similar way as Proposition \ref{inter-CM}.
\end{proof}

There is a decomposition \[\Xi = \bigsqcup_{z \in W_{S, T, r}/W_{M, r}} \Xi_z,\] where \begin{align*}\Xi_z &= \{(x, \psi) \in \Xi; \psi \in \CI_r z P_r / V_r\} \\ &=\{(x, vzm V_r) \in \CI_{\phi, r}^+ \times \CI_{z, r} z P_r/V_r; F(x v z m V_r) = v z m V_r\} \\ &\cong \{(x, v, zm V_r) \in \CI_{\phi, r}^+ \times \CI_{z, r} \times zP_r/V_r; F(x v z m V_r) = v z m V_r\}. \end{align*}

Let $\psi \in \CI_r z P_r/V_r$. Consider the fiber of the projection $\pr_2: \Xi_z \to G_r/V_r$ to the second factor:  \[\Xi_z(\psi) = \pr_2\i(\psi) = \{x  \in \CI_{\phi, r}^+; F(x \psi) = \psi\}.\]

\begin{lemma} \label{Xi-deform}
    We have $R\G_c(\Xi_z(\psi), \CF) = 0$ for $\psi \in \CI_r z P_r/V_r \sm \CI_{\phi, r}^+ z P_r/V_r$.
\end{lemma}
\begin{proof}
    Write $\psi = v z m V_r$ with $v \in \CI_{z, r}$ and $m \in M_r$. Let \[h = \max\{f^\dag; f < s_f, v_f \neq 0\} \in E_\phi,\] which exists by the assumption on $\psi$. Moreover, $\Im u_h \subseteq \CI_r \cap {}^z V_r$.
    
    Let $\d: \BA^1 \to G_r$ given by $\t \mapsto v u_h(v_{h^\dag}\i \t) v\i$. By Lemma \ref{commutator}, we have \[\d \equiv \l_{\a_h^\vee, r_h} \mod \hCI_{\phi, r}^{\dag}.\] In particular, $\Im \d \subseteq \CI_{\phi, r}^+$. Thus, there is an isomorphism \[\D: \BA^1 \times \Xi_z(\psi) \to \BA^1 \times \Xi_z(\psi), \quad (\t, x) \mapsto (\t, x \d(\t)).\] Then we have \[\D^*(\ov\BQ_\ell \boxtimes \CF) \cong (\l_{\a_h^\vee, r_h}^*\CF) \boxtimes \CF.\] By a similar argument in the proof of Lemma \ref{deform}, we deduce that $R\G_c(\Xi_z(\psi), \CF) = 0$ as desired.
\end{proof}

Let \[\Xi_z^0 = \pr_2\i(\CI_{\phi, r}^+ z P_r/V_r) \cong \{(x, v, zm V_r) \in \CI_{\phi, r}^+ \times (\CI_{\phi, r}^+ \cap {}^z \ov V_r) \times zP_r/V_r; F(x z m V_r) = vzm V_r\},\] on which $T_0^F$ acts by \[t: (x, v, zm V_r) \to (t\i x t, t\i v t, t\i zm V_r).\]
\begin{corollary} \label{Xi-zero}
    We have $R\G_c(\Xi_z^0, \pr_1^*\CF)[\phi_0] \cong R\G_c(\Xi_z, \pr_1^*\CF)[\phi_0]$.
\end{corollary}
\begin{proof}
    It follows from Lemma \ref{Xi-deform}.
\end{proof}

\subsection{The vanishing result on $\Xi_z^0$}
Consider the cover of $\Xi_z^0$ defined by \[\hXi_z^0 = \{(x, \psi, s) \in \Xi_z^0 \times T_{0+:r}; \t_r(x) = L(s)\},\] on which $T_r^F$ acts by \[(t, \t): (x, \psi, s) \mapsto (t\i x t, t\i \psi, s\t\i),\] where $(t, \t) \in T_0^F \times T_{0+:r}^F \cong T_r^F$.

\begin{lemma} \label{Xi-cover}
We have $R\G_c(\Xi_z^0, \pr_1^*\CF)[\phi_0] \cong R\G_c(\hXi_z^0, \ov\BQ_\ell)[\phi]$.    
\end{lemma}

For a cocharacter $\chi$ of $G$ we define $N_\chi: \BG_m \to G_r$ by $\t \mapsto N_F^{F^N} \circ \chi([\t])$. Recall that $N \in \BZ_{\ge 1}$ such that $F^N$ acts trivially on $X_*(T)$. Let $\mu = z(\l)$.
\begin{lemma} \label{empty}
    There is an action of $\BG_m$ on $\hXi_z^0$ given by \[\t: (x, v, zmV_r, s) \mapsto ({}^{N_\mu(\t)} x , {}^{F(N_\mu(\t))} v, F(N_\mu(\t)) zm F(N_\l(\t\i)) V_r, s).\] As a consequence, $(\Xi_z^0)^{\BG_m} = \emptyset$ and $H^*_c(\hXi_z^0, \ov\BQ_\ell)[\phi\i] = 0$.
\end{lemma}
\begin{proof}
    Note that $\t_r({}^{N_\mu(\t)} x) = \t_r(x)$. By \cite[Lemma 7.3]{LNY26}, the action of $\BG_m$ on $\hXi_z^0$ is well-defined, and the remaining statement follows in the same way.
\end{proof}

\begin{theorem}
    Assume $B = T U$ is convex. Then $(R_{T_r}^{G_r}(\phi))^{V_r^F} = 0$. 
\end{theorem}
\begin{proof}
    It follows from Lemma \ref{V-geo}, Corollary \ref{Xi-zero}, Lemma \ref{Xi-cover} and Lemma \ref{empty}.
\end{proof}

\end{document}